\documentclass[12pt,a4paper,reqno]{amsart}
\usepackage{amsmath}
\usepackage{amsfonts}
\usepackage{amssymb}
\numberwithin{equation}{section}

\theoremstyle{plain}

\newtheorem{theorem}[subsection]{Theorem}
\newtheorem{proposition}[subsection]{Proposition}
\newtheorem{lemma}[subsection]{Lemma}
\newtheorem{corollary}[subsection]{Corollary}

\theoremstyle{definition}

\newtheorem{definition}[subsection]{Definition}

\newtheorem{example}[subsection]{Example}
\newtheorem{remark}[subsection]{Remark}
\newtheorem{remarks}[subsection]{Remarks}

\renewcommand{\leq}{\leqslant}
\renewcommand{\geq}{\geqslant}

\newsavebox{\proofbox}
\savebox{\proofbox}{\begin{picture}(7,7)%
  \put(0,0){\framebox(7,7){}}\end{picture}}

\def\proof{\noindent\textit{Proof. }}
\def\endproof{\hfill{\usebox{\proofbox}}}

\def\N{\mathbb{N}}

\def\Z{\mathbb{Z}}
\def\R{\mathbb{R}}

\def\vol{\operatorname{vol}}

\def\size{\mbox{size}}
\def\sizelem{\mbox{\emph{size}}}
\def\dfrei{d_{\mbox{\scriptsize Frei}}}
\def\dfreilem{d_{\mbox{\scriptsize \emph{Frei}}}}

\def\ni{\noindent}
\def\vs{\vspace{11pt}}

\begin{document}

\title[Freiman--Bilu]{Compressions, Convex Geometry and the Freiman-Bilu Theorem}

\author{Ben Green}
\address{Department of Mathematics, University of Bristol, University Walk, Bristol BS8 1TW, England}
\email{b.j.green@bristol.ac.uk}

\author{Terence Tao}
\address{Department of Mathematics, UCLA, Los Angeles CA 90095-1555}
\email{ tao@math.ucla.edu}

\thanks{The first author is a Clay Research Fellow, and is pleased to acknowledge the support of the Clay Mathematics Institute. Some of this work was carried out while he was on a long-term visit to MIT, and he would like to express his gratitude for the institute's hospitality.  The second author is supported by a grant from the Packard Foundation.}

\begin{abstract} 
\ni We note a link between combinatorial results of Bollob\'as and Leader concerning sumsets in the grid, the Brunn-Minkowski theorem and a result of Freiman and Bilu concerning the structure of sets $A \subseteq \Z$ with small doubling.\vs

\ni Our main result is the following. If $\varepsilon > 0$ and if $A$ is a finite nonempty subset of a torsion-free abelian group with $|A + A| \leq K|A|$, then $A$ may be covered by $e^{K^{O(1)}}$ progressions of dimension $\lfloor \log_2 K + \varepsilon \rfloor$ and size at most $|A|$.\vs

\end{abstract}

\maketitle

\section{Introduction}

\ni A famous theorem of Freiman \cite{freiman} concerns the structure of finite non-empty sets $A \subseteq \Z^m$ with small doubling property, that is to say sets for which the \emph{doubling constant}
\[ \sigma[A] := \frac{|A + A|}{|A|}\] is small in comparison with the cardinality $|A|$. Here of course $A+A := \{ a_1+a_2: a_1,a_2 \in A \}$ denotes the Minkowski sum of $A$ with itself.  Freiman's theorem asserts that such sets must be 
contained in a \emph{generalised arithmetic progression} (or simply \emph{progression}, for short) of small dimension and size.

\begin{definition}[Progressions]
Let $d$ and $L_1,\dots,L_d$ be positive integers. Then a \emph{progression} of dimension $d$ is a set of the form \[
 P := x_0 + [L_1] \cdot x_1 + \ldots + [L_d] \cdot x_d = \{ x_0 + \sum_{i=1}^d \mu_i x_i  : \mu_i \in \Z, \mu_i \in [L_i]\},\]
where we write $[n] := \{0,\dots,n-1\}$.
The \emph{size} (or \emph{volume}) of $P$ is defined by
\[ \size(P) := L_1 \dots L_d.\]
If all the sums $x_0 + \sum_{i=1}^d \mu_i x_i$ in $P$ are distinct then we say that $P$ is \emph{proper}, in which case $\size(P) = |P|$. More generally if $t \geq 1$ is an integer, and if all the sums in 
\[ tP := \{x_0 + \sum_{i=1}^d \mu_i x_i  : \mu_i \in \Z, \mu_i \in [tL_i]\}\]
are distinct then we say that $P$ is $t$-proper.\hfill$\diamond$
\end{definition}

\ni Freiman proved that there are functions $d$ and $F$ such that any set $A \subseteq \Z^m$ is contained in a progression of dimension at most $d(\sigma[A])$ and size at most $F(\sigma[A])|A|$. Chang \cite{chang-freiman}, building on earlier ideas of Ruzsa \cite{ruzsa-freiman}, proved a rather effective version of Freiman's theorem in which the functions $d$ and $F$ are close to being best possible. Before giving our statement of Chang's theorem, we need to recall some definitions concerning Freiman homomorphisms.

\begin{definition}[Freiman homomorphisms and dimension]
Let $G,G'$ be abelian (additive) groups, and suppose that $A \subseteq G$ and $A' \subseteq G'$ are sets. We say that a map $\phi : A \rightarrow A'$ is a \emph{Freiman homomorphism} if whenever $a_1,a_2,a_3,a_4 \in A$ are such that 
\[ a_1 + a_2 = a_3 + a_4\]
we have 
\[ \phi(a_1) + \phi(a_2) = \phi(a_3) + \phi(a_4).\]
If $\phi$ has an inverse which is also a Freiman homomorphism, then we say that $\phi$ is a \emph{Freiman isomorphism}. 
We define the \emph{Freiman dimension} $\dfrei(A)$ of $A$ to be the largest dimension $d$ for which $A$ is Freiman isomorphic to a subset of $\Z^d$ not contained in a proper affine subspace. 
\end{definition}

\ni The following version of Freiman's theorem was essentially proved by Chang. See also \cite{green-edinburgh-mit}.

\begin{proposition}[Chang]\label{chang-thm}\cite{chang-freiman}
Let $t \geq 1$ be an integer, and let $A$ be a finite subset of a torsion-free abelian group \textup{(}such as $\Z^m$\textup{)} with $|A| \geq 2$ and $\sigma[A] = K$. Then there is a $t$-proper progression $P$ with \[ \dim(P) \leq \dfreilem(A)\] and 
\[ |P| = \sizelem(P) \leq t^K\exp(CK^2 \log^3 K)|A|.\]
\end{proposition}

\begin{remarks} In this paper $t$ will always be taken to be $2$, in which case the $t^K$ can be absorbed into the $\exp(CK^2 \log^3 K)$ term. Here, $C$ denotes an absolute constant which may change from line to line.  Theorem \ref{chang-thm} is stated in an arbitrary torsion-free abelian group, but as $A$ is finite one can easily reduce to the case when the ambient group is 
finitely generated and hence isomorphic to $\Z^m$ for some $m$.  In fact by the machinery of Freiman isomorphisms one can reduce further, to the case when $A$ is a subset of the integers.  The case of groups with torsion has some further subtleties; see \cite{gr-4}.
\end{remarks}

\ni As stated, Proposition \ref{chang-thm} does not quite imply Freiman's theorem, since we have not bounded $\dfrei(A)$ in terms of $\sigma[A]$. To achieve this we need a further (elementary) lemma of Freiman:

\begin{proposition}[Freiman's lemma]\label{frei-lem}\cite{freiman}
Let $A \subseteq \Z^m$ be a finite non-empty set, and write $d := \dfreilem(A)$. We have the inequality 
\[ |A + A| \geq (d + 1)|A| - \textstyle \frac{1}{2} \displaystyle d(d+1).\]
As a consequence, if $K = \sigma[A]$, $\varepsilon > 0$ and $|A| \geq CK^2/\varepsilon$ then we have the bound \begin{equation}\label{strong-dim-bound} d \leq \lfloor K - 1 + \varepsilon\rfloor.\end{equation}
\end{proposition}

\ni See also \cite{ruzsa-multid} or \cite[Section 5.2]{tao-vu} 
for some further discussion of this result and a generalization to sums $A+B$ of different sets.

\begin{example}\label{changex} Chang's result is fairly close to best possible. To see why, consider a set $A$ of the form
\begin{equation}\label{eq1} A := \bigcup_{i = 1}^K (x_i + \{1,\dots,m\}),\end{equation}
where the $x_i$ form a very lacunary sequence of integers (or are linearly independent in $\R^m$). A moment's thought confirms that if $m \ggg K$ then $|A + A| \approx (K+1)|A|$ and $\dfrei(A) = K$. Yet $A$ is not contained in any arithmetic progression of dimension smaller than $K$, and if $P$ is a progression of dimension $K$ containing $A$ then $|P| \geq \frac{2^{K-1}}{K}|A|$. 
\end{example}

\ni One distressing feature of this example is that it places a limit on how close Freiman's theorem can be to a necessary and sufficient condition for small doubling. It is not hard to see that Chang's theorem gives a double implication

\begin{equation}\label{eq2} \mbox{($\sigma[A] \leq K$)} \;\; \Longrightarrow \;\; \mbox{($A$ has structure)} \;\; \Longrightarrow \;\; \mbox{($\sigma[A] \leq f(K)$)}.\end{equation}

\ni However, the function $f(K)$ is of order $\exp(K^{2 + \epsilon})$, which is much larger than $K$. Example 1 shows that such a poor dependence is more-or-less necessary if ``structure'' means ``is contained in a progression''. It would be nice to find a notion of ``structure'' which allowed for an implication of the type \eqref{eq2} in which $f(K)$ was merely polynomial. \vs

\ni It is quite possible that there \emph{is} such a variant of Freiman's theorem. One approach might be to cover $A$ by a few progressions of dimension merely $O(\log K)$. Ruzsa \cite{ruzsa-fin-fields} attributes the analogue of this idea in the finite field setting to Katalin Marton. The notion that one can find an efficient covering of this type is known as the Polynomial Freiman-Ruzsa Conjecture (PFR). It has only been properly formulated in the finite field setting (see \cite{green-fin-field}), and indeed some care is needed when working in $\Z$. In a future paper we will give an example to show that covering by progressions is not a general enough notion of structure to allow $f(K)$ to be polynomial. One needs to work instead with unions of \emph{convex progressions}, which are essentially the image of the set of lattice points in a convex body under an affine map.\vs

\ni In this direction we have a remarkable result of Freiman and Bilu \cite{bilu,bilu2}, which shows that
given any finite subset $A$ of a torsion-free abelian group with $\sigma[A]=K$, and given any $\varepsilon > 0$, one can
find a progression $P'$ with dimension at most $\lfloor \log_2 K + \varepsilon \rfloor$ and size at most $|A|$ such that
$$ |A \cap P'| \geq \exp(-\exp(C_\varepsilon K^{C_\varepsilon})) |A|,$$
where the constants $C_\varepsilon > 0$ depend on $\varepsilon$.
Thus $A$ has large intersection with a low-dimensional progression.
Our main result is to improve the bounds in this result, and to show that $A$ is in fact covered by relatively few translates
of a low-dimensional progression\footnote{There is actually not much difference between the statements that $A$ has large intersection with a progression, and that $A$ can be covered by few translates of a progression, thanks to sum set estimates and 
the Ruzsa covering lemma; see \cite{ruzsa-freiman} or \cite[Chapter 2]{tao-vu}.}.  More precisely, we have

\begin{theorem}\label{mainthm}  Let $A$ be a finite subset of a torsion-free abelian group with $\sigma[A] = K$ and $|A| \geq 2$. Then for
any $\varepsilon \in (0,1]$ one can cover $A$ by at most $\exp(C K^3 \log^3 K) / \varepsilon^{CK}$ 
translates of a progression $P'$ with dimension at most $\lfloor \log_2 K + \varepsilon \rfloor$ and size at most $|A|$.
\end{theorem}

\begin{remark} For some purposes, such as the papers of Szemer\'edi and Vu \cite{szemeredi-vu1,szemeredi-vu2,szemeredi-vu3}, it is necessary to have the Freiman-Bilu theorem for rather small $\varepsilon$ (say $\varepsilon < 1/2$). However for many applications in additive combinatorics
(for example to the problem of obtaining inverse theorems for the Gowers $U^3$-norm \cite{gowers,green-tao-gowersu3} or to the sum-product problem \cite{bourgain-chang}) we would be just as happy with $\varepsilon = 1$, or even to permit $P'$ to have dimension as large as $100 \log K$ (say).  
\end{remark}

\begin{remark} In this remark we set $\varepsilon=1$ and assume that $K$ is large for 
simplicity. It is of interest to understand 
whether the quantity $\exp( C K^3 \log^3 K)$ in the above theorem can be improved significantly.  We do not expect to be able
to obtain a polynomial bound such as $K^C$, unless one generalizes the notion of progression to also include the convex progressions mentioned earlier.  A weaker bound such as $K^{C (\log\log K)^C}$ might however be possible (though probably quite difficult).
We will return to these issues in a future paper. 
\end{remark}

\begin{remark} One can ensure that the progression $P'$ appearing in Theorem \ref{mainthm} is proper by standard rank reduction
methods: see, for example, \cite{chang-freiman} or \cite[Section 3.6]{tao-vu}.  The dimension $d \leq \lfloor \log_2 K + \varepsilon \rfloor$ 
of $P'$ is so small that the losses incurred by such a rank reduction procedure (which are generally of the form $\exp(d^C)$ or so)
will be negligible when compared against the factor of $\exp( C K^3 \log^3 K)$ in the theorem; we omit the details.  
In fact one can even make $P'$ $t$-proper for a very large value of $t$ (such as $t \sim \exp(C K^3)$) without any significant degradation of constants.  However, one would need to be significantly more careful with these issues if one wanted to obtain polynomial type bounds rather than exponential ones.
\end{remark}

\section{A preliminary doubling lemma}

\ni The purpose of this section is to prove a variant of Freiman's lemma (Proposition \ref{frei-lem}), which essentially 
gives an exponential doubling constant $2^d$, but requires that the set $A$ is a large subset of a $d$-dimensional box.
More precisely, we shall prove

\begin{theorem}[Doubling in boxes]\label{box-cor}
Suppose that $A \subseteq [L_1] \times \ldots \times [L_d]$. Then
\[ |A + A| \geq 2^d |A| + \prod_{i=1}^d (2L_i-1) - \prod_{i=1}^d 2L_i.\]
\end{theorem}

\begin{example} From this theorem and the mean value theorem we see that if 
$A \subseteq [k]^d$, then $|A+A| \geq 2^d |A| - d (2k)^{d-1}$.  Thus if $|A| \gg d(2k)^{d-1}$ then $A$ will have
doubling constant close to or greater than $2^d$.
\end{example}

\ni The proof of Theorem \ref{box-cor} uses a simple case of the Brunn-Minkowski inequality.

\begin{proposition}[Brunn-Minkowski inequality]
Suppose that $A$ and $B$ are bounded open sets in $\R^d$. Then
\[ \vol(A + B) \geq (\vol(A)^{1/d} + \vol(B)^{1/d})^d.\]
In particular we have
$$ \vol(A+B) \geq 2^d \min( \vol(A), \vol(B) ).$$
\end{proposition}

\ni See \cite{brunn-minkowski-survey} for an excellent survey of this inequality
and related ideas. 
The following is a discrete consequence of the Brunn-Minkowski inequality. 

\begin{lemma}[Discretized Brunn-Minkowski]\label{disc-bm}
Suppose that $X,Y \subseteq \R^m$ are finite sets.  Then for any $d \leq m$ we have
\[ |X + Y + \{0,1\}^d| \geq 2^d \min(|X|,|Y|),\]
where we embed $\Z^d$ \textup{(}and hence $\{0,1\}^d$\textup{)} inside $\R^m$ in the obvious manner.
\end{lemma}

\proof 
Let us first prove the claim when $X,Y$ lie in $\Z^d$.  Set $A := X + (0,1)^d$ and
$B := Y + (0,1)^d$. Then $A,B$ are bounded open sets with $A+B = X+Y+\{0,1\}^d+(0,1)^d$.  From the Brunn-Minkowski inequality
and the observation that $\vol(X + (0,1)^d) = |X|$ whenever $X$ is a subset of $\Z^d$, we obtain
$$ |X + Y + \{0,1\}^d| \geq 2^d \min(|X|,|Y|)$$
as desired.\vs

\ni By translation we see that the same argument applies when $X$ and $Y$ live in (possibly distinct) cosets of $\Z^d$.  Now we consider
the general case when $X,Y$ are arbitrary finite subsets of $\R^m$.  We split $X = X_1 \cup \ldots \cup X_k$ and $Y = Y_1 \cup \ldots \cup Y_l$, where $X_1,\ldots,X_k$ live in disjoint cosets of $\Z^d$ and similarly for $Y_1,\ldots,Y_l$.  Without loss of
generality we may assume that $|X_1| \geq |X_i|, |Y_j|$ for all $i$, $1 \leq i \leq k$, and $j$,  $1 \leq j \leq l$.  From the preceding arguments we have \[ |X_1 + Y_j + \{0,1\}^d| \geq 2^d |Y_j| \] for all $j$, $1 \leq j \leq l$.  Since the sets $X_1 + Y_j + \{0,1\}^d$ lie in disjoint cosets of $\Z^d$ as $j$ varies, we conclude upon taking unions that \[ |X+Y+\{0,1\}^d| \geq 2^d |Y| \geq 2^d \min(|X|,|Y|), \]
as desired.
\endproof\vs

\ni The lemma is sharp, as can be seen by setting $X=Y$ equal to a box such as $[k]^d$.  However, we would still 
like to eliminate the $\{0,1\}^d$ term on the left-hand side.  For some refined estimates concerning this problem (especially in the case when $X$ and $Y$ have very different sizes), see \cite{gardner-gronchi}.  However, following Bollob\'as and Leader \cite{bollobas-leader}, we 
shall pursue a simpler approach, using the machinery of down-sets and compressions. The idea of using compressions to study sumsets seems to have first appeared in the work of Freiman \cite{freiman}. In particular the theorem on page 27 of that book is a special case of our Lemma \ref{doub}.

\begin{definition}[Down-sets]  We use $\N_0 = \{0,1,\ldots\}$ to denote the non-negative integers.
A subset $B \subseteq \N_0^d$ is said to be an \emph{$i$-down-set} for some $1 \leq i \leq d$
if whenever $(x_1,\ldots,x_d) \in B$, then $(x_1,\ldots,x_{i-1},y_i,x_{i+1},\ldots,x_d) \in B$ for all
integers $0 \leq y_i \leq x_i$.  We say that $B$ is a \emph{down-set} if it is an $i$-down-set for each $1 \leq i \leq d$.
\end{definition}

\begin{example} The cube $[k]^d$ is a down-set, and more generally the box $[L_1] \times \ldots \times [L_d]$ is a down-set.
If $A$ and $B$ are down-sets, then $A+B$ is also a down-set.
\end{example}

\ni We can remove the $\{0,1\}^d$ term in Lemma \ref{disc-bm} for down-sets, at the cost of a small error term.  For each
$I \subseteq \{1,\ldots,d\}$ let $\pi_I: \R^d \to \R^d$ be the orthogonal projection $\pi_I(\sum_{i=1}^d x_i e_i) := \sum_{i \in I} x_i e_i$, where $e_1,\ldots,e_d$ is the standard basis of $\R^d$.  Let us observe the identity
\begin{equation}\label{xoi}
|X + \{0,1\}^d| = |X - \{0,1\}^d| = \sum_{I \subseteq  \{1,\ldots,d\}} |\pi_I(X)|
\end{equation}
whenever $X$ is a down-set.  The first identity follows since $-\{0,1\}^d$ is a translate of $\{0,1\}^d$.  The second identity follows from the easily verified observation that the sets 
$$\pi_I(X) - \sum_{i \in \{1,\ldots,d\} \backslash I} e_i = \{ (x_1,\ldots,x_d) \in X - \{0,1\}^d: x_i = -1 \hbox{ if and only if } i \not \in I \}$$
partition $X - \{0,1\}^d$.\vs

\ni From Lemma \ref{disc-bm} and \eqref{xoi} applied to the set $X+Y$, we conclude

\begin{corollary}\label{coco} Let $X,Y \subseteq \N_0^d$ be down-sets.  Then
$$ |X + Y| \geq 2^d \min(|X|,|Y|) - \sum_{I \subsetneq \{1,\ldots,d\}} |\pi_I(X+Y)|.$$
\end{corollary}

\ni We now use the method of compressions to generalize the above lemma to \emph{subsets} of down-sets, rather than down-sets themselves.

\begin{lemma}\label{doub}  Let $X, Y \subseteq \N_0^d$ be down-sets, and let $A \subseteq X$, $B \subseteq Y$ be arbitrary \textup{(}and so not necessarily down-sets\textup{)}.  Then
\begin{equation}\label{ab}
 |A + B| \geq 2^d \min(|A|,|B|) - \sum_{I \subsetneq \{1,\ldots,d\}} |\pi_I(X+Y)|.
 \end{equation}
\end{lemma}

\proof For each $1 \leq i \leq d$, define the \emph{$i$-compression} ${\mathcal C}_i(A)$ of $A$ to be the set formed by moving
the elements of $A$ in the $-e_i$ direction until one obtains an $i$-down-set.  More precisely, one has
$$ {\mathcal C}_i(A) := \{ (x_1,\ldots,x_{i-1},t,x_{i+1},\ldots,x_d): 0 \leq t < |A \cap \pi_{\{i\}}^{-1}(x_1,\ldots,x_{i-1},0,x_{i+1},\ldots,x_d)| \}.$$
For instance, if $A = \{(0,1), (1,0), (2,0) \}$, then ${\mathcal C}_1(A) = \{ (0,1), (0,0), (1,0) \}$
and ${\mathcal C}_2(A) = \{ (0,0), (1,0), (2,0) \}$.  Note that this example also shows that ${\mathcal C}_i$ and ${\mathcal C}_j$
do not in general commute with each other.\vs

\ni One can easily verify that if $A \subseteq X$, then ${\mathcal C}_i(A) \subseteq X$ has the same cardinality as $A$, 
and that ${\mathcal C}_i(A)$ is an $i$-down-set. Furthermore, if $A$ was already a $j$-down-set for some $j \neq i$,
then ${\mathcal C}_i(A)$ is still a $j$-down-set.  Finally, we observe the crucial compression property
$$ |{\mathcal C}_i(A) + {\mathcal C}_i(B)| \leq |A+B|.$$
To verify this, it suffices to show that
$$ |{\mathcal C}_i(A) + {\mathcal C}_i(B) \cap \pi_{\{i\}}^{-1}(x)| \leq |(A+B) \cap \pi_{\{i\}}^{-1}(x)|$$
for all $x \in \pi_{\{i\}}(\N_0^d)$.  But the right-hand side can be written as
$$ | \bigcup_{y+z = x} (A \cap \pi_{\{i\}}^{-1}(y)) + (B \cap \pi_{\{i\}}^{-1}(z))|$$
while the left-hand side can be computed as
$$ \max( \sup_{y+z=x} |A \cap \pi_{\{i\}}^{-1}(y)| + |B \cap \pi_{\{i\}}^{-1}(z)| - 1, 0),$$
and the claim follows from the easy estimate $|X+Y| \geq \max( |X|+|Y|-1, 0 )$.\vs

\ni From the above discussion we see that to verify the estimate \eqref{ab} it suffices to do so with $A$ and $B$ replaced
by ${\mathcal C}_i(A)$ and ${\mathcal C}_i(B)$ respectively.  Iterating this observation, we may replace $A$ and $B$ by the ${\mathcal C}_d \ldots {\mathcal C}_1(A)$ and ${\mathcal C}_d \ldots {\mathcal C}_1(B)$ respectively.  But these latter sets are down-sets, and the claim now follows from Corollary \ref{coco}.
\endproof\vs

\ni We can now prove Theorem \ref{box-cor}.  We apply Lemma \ref{doub} with $A=B$ and $X=Y=[L_1] \times \ldots \times [L_d]$, noting that $X+Y = [2L_1-1] \times \ldots \times [2L_d-1]$.  The claim then follows after using \eqref{xoi} to compute
$\sum_{I \subsetneq \{1,\ldots,d\}} |\pi_I(X+Y)|$.
\endproof

\section{Proof of Theorem \ref{mainthm}}

\ni Let $A$ be such that $\sigma[A] = K$. Take $t = 2$ in Chang's theorem, and suppose that 
\[ P := \{x_0 + \mu_1 x_1 + \dots + \mu_d x_d : \mu_i \in \Z, \mu_i \in [L_i]\}\] is the resultant $2$-proper progression containing $A$. We may assume that $|A| \geq CK^2/\varepsilon$ since the claim is trivial otherwise (covering $A$ by translates of a point).
Chang's theorem then allows us to take 
\begin{equation}\label{d-bound}d = \lfloor K - 1 + \varepsilon\rfloor.
\end{equation}
Actually we will not need such a precise bound; $d \leq 2K$ will be more than sufficient. Observe that the map $\phi : P \rightarrow [L_1] \times \ldots \times [L_d]$ defined by 
\[ \phi\big(x_0 + \sum_{i=1}^d \mu_i x_i \big) = (\mu_1,\dots,\mu_d)\] is a Freiman isomorphism. In particular we have $|A + A| = |\phi(A) + \phi(A)|$ and $|A| = |\phi(A)|$.  Note that Chang's theorem also gives
\begin{equation}\label{a-lower}
|A| \geq \exp(-CK^2 \log^3 K) L_1 \dots L_d.
\end{equation}

\ni Without loss of generality we may order the $L_i$ so that $L_1 \geq L_2 \geq \dots \geq L_d$. Fix $\varepsilon > 0$, and write 
\[ l := \lfloor \log_2 K + \varepsilon \rfloor.\] At the heart of our argument is a proof that $L_{l+1}$ (and hence each of $L_{l+2},\dots,L_d$) is small. For each $x\in [L_{l+2}] \times \dots \times [L_d]$, define the fibre $\phi(A)_x$ of $\phi(A)$ by
\[ \phi(A)_x := \phi(A) \cap \big( [L_1] \times \dots \times [L_{l+1}] \times \{x\}\big).\]
The sets $\phi(A)_x + \phi(A)_x$ are disjoint, and we thus have
\[ |A+A| = |\phi(A) + \phi(A)| \geq \sum_x |\phi(A)_x + \phi(A)_x|. \]
Now Theorem \ref{box-cor} easily implies that
\[ |\phi(A)_x + \phi(A)_x| \geq 2^{l+1}|\phi(A)_x| - d 2^d L_1 \dots L_{l} \] and hence, summing over $x$, we have
\begin{equation}\label{eq11} |A + A| \geq 2^{l+1} |A| - d 2^d L_1 \dots L_{l} L_{l+2} \dots L_d.\end{equation}
However we know that $|A + A| \leq K|A|$, and from the definition of $l$ we clearly have $2^{l+1} \geq 2^{\varepsilon}K \geq  (1 + \varepsilon/2)K$. Substituting into \eqref{eq11} we obtain
\[ \frac{\varepsilon |A|}{2} \leq d 2^d L_1 \dots L_{l} L_{l+2} \dots L_d;\]
combining this with \eqref{d-bound} and \eqref{a-lower} yields
\begin{equation}\label{eq12} L_{l+1} \leq \exp(C K^2 \log^3 K)/\varepsilon.\end{equation}
From this, the monotonicity of the $L_j$, and \eqref{d-bound} we see that we can cover the box $[L_1] \times \ldots \times [L_d]$
(and hence $\phi(A)$) by at most $\exp(C K^3 \log^3 K) / \varepsilon^{CK}$ copies of $[L_1] \times \ldots \times [L_l]$.  Applying the inverse map $\phi^{-1}$,
we can cover $A$ by at most $\exp(C K^3 \log^3 K) / \varepsilon^{CK}$ translates of a progression $\tilde P$ of dimension $l$ and volume 
$L_1 \ldots L_l$. Now from \eqref{d-bound}, \eqref{a-lower} and \eqref{eq12} we have $|A| \geq \varepsilon^{-CK}\exp(C K^3 \log^3 K) L_1 \ldots L_l$, so we can easily cover
$\tilde P$ by at most $\exp(C K^3 \log^3 K) / \varepsilon^{CK}$ translates of a progression $P$ of dimension $l$ and volume at most $|A|$. The claim follows.
\endproof

\section{Exponential doubling and cubes}

\ni In this section we present another variant on the above results, which asserts that any set which contains a non-degenerate
$d$-dimensional parallelopiped will have doubling constant exponential in $d$ or larger.

\begin{proposition}\label{exp-cor} Suppose that $A \subseteq  \R^m$ is a finite set which contains a non-degenerate parallelopiped 
\[ P := v_0 + \{0,1\} \cdot v_1 + \dots + \{0,1\} \cdot v_d,\]
for some linearly independent $v_1,\dots,v_d \in \R^m$. Then $|A + A| \geq 2^{d/2}|A|$.
\end{proposition}

\proof By applying a suitable linear transformation and translation we may assume that $P = \{0,1\}^d \subseteq \Z^d \subseteq \R^m$. 
We have $|A+A| = \sigma[A] |A|$.  By the Pl\"unnecke inequality for commutative graphs \cite{plun} (see also \cite{nath, ruzsa-plun, tao-vu}) there thus exists a non-empty subset $B$ of $A$ such that $|B+A+A| \leq \sigma[A]^2 |B|$.  On the other hand, 
from Lemma \ref{disc-bm} we have
$$ |B+A+A| \geq |B+A+\{0,1\}^d| \geq 2^d |B|.$$
Hence $\sigma[A] \geq 2^{d/2}$, as desired.
\endproof\vs

\begin{remark} This bound is somewhat close to being sharp; if $A = \{0,1\}^d$ then $|A+A| = (3/2)^d |A|$.  By considering sets
such as $A = \{0,1\}^{d-1} \times [k]$ for large $k$ we see that similar behaviour persists even if $A$ is allowed to be arbitrarily large.  Thus we cannot hope to force a doubling constant close to $2^d$ from the hypothesis that $A$ contains a $d$-dimensional parallelopiped, in contrast to Theorem \ref{box-cor}.  As a kind of converse to the above result, we have the Freiman cube lemma (see e.g. \cite{bilu} or \cite[Corollary 5.19]{tao-vu}), which asserts that if $A \subseteq \R^d$ has doubling constant $\sigma[A] = K$, then $A$ contains a $d$-dimensional parallelopiped which captures a significant fraction (such as $(3K)^{-2^d}$) of $A$ in its convex hull; however this parallelopiped may be degenerate.
\end{remark}

\ni It is somewhat diverting to compare this with Freiman's lemma (Proposition \ref{frei-lem}), which may be stated in the following form.

\begin{proposition}[Freiman's lemma, again]
Suppose that $A \subseteq  \R^m$ is a finite set which contains some non-degenerate simplex 
\[ P := \{v_0 , v_0 + v_1,\dots, v_0 + v_d\}\]
where $v_1,\dots,v_d$ are linearly independent. Then $|A + A| \geq (d+1)|A| - \frac{1}{2}d(d+1)$.
\end{proposition}

\ni The maximal possible value of $d$ in Freiman's lemma (applying Freiman isomorphisms if necessary) is precisely $\dfrei(A)$. Because of this connection, it is tempting to think of the maximum permissible $d$ in Corollary \ref{exp-cor} as the \emph{exponential Freiman dimension} of $A$. We do not, at this stage, envisage any applications of it.

\end{document}